\documentclass[onecolumn,english,letterpaper,11pt]{amsart}
\usepackage[latin9]{inputenc}
\usepackage{amsthm}
\PassOptionsToPackage{normalem}{ulem}
\usepackage{ulem}

\numberwithin{equation}{section} 
\numberwithin{figure}{section} 
\theoremstyle{plain}

\usepackage{fullpage}
\usepackage{babel}
\usepackage{amsmath}
\usepackage{url}
\usepackage{color}
\usepackage{graphics}
\usepackage{graphicx}
\usepackage{algorithm2e}
\usepackage[version=3]{mhchem}
\usepackage{lineno}

\usepackage[bitstream-charter]{mathdesign}

\newtheorem{theorem}{Theorem}[section]
\newtheorem{lemma}[theorem]{Lemma}

\newtheorem{definition}[theorem]{Definition}

\newtheorem{example}[theorem]{Example}

\newtheorem{remark}[theorem]{Remark}

\newcommand{\Z}{{\mathbb Z}}

\newcommand{\HPS}{\mathcal H}

\newcommand{\1}{{\bf 1}}

\DeclareMathOperator{\supp}{supp}
\DeclareMathOperator{\lt}{LT}
\DeclareMathOperator{\NF}{NF}
\DeclareMathOperator{\cluster}{{\mathcal C}}
\DeclareMathOperator{\clustergraph}{{\mathcal G}}

\def\ve#1{\mathchoice{\mbox{\boldmath$\displaystyle\bf#1$}}
{\mbox{\boldmath$\textstyle\bf#1$}}
{\mbox{\boldmath$\scriptstyle\bf#1$}}
{\mbox{\boldmath$\scriptscriptstyle\bf#1$}}}

\usepackage{enumerate}
\def\ve#1{\mathchoice{\mbox{\boldmath$\displaystyle\bf#1$}}
{\mbox{\boldmath$\textstyle\bf#1$}}
{\mbox{\boldmath$\scriptstyle\bf#1$}}
{\mbox{\boldmath$\scriptscriptstyle\bf#1$}}}
\newcommand\vea{{\ve a}}
\newcommand\veb{{\ve b}}

\newcommand\ved{{\ve d}}
\newcommand\vece{{\ve e}}

\newcommand\vem{{\ve m}}
\newcommand\vep{{\ve p}}

\newcommand\ver{{\ve r}}
\newcommand\ves{{\ve s}}
\newcommand\vet{{\ve t}}
\newcommand\veu{{\ve u}}
\newcommand\vev{{\ve v}}
\newcommand\vew{{\ve w}}
\newcommand\vex{{\ve x}}
\newcommand\vey{{\ve y}}
\newcommand\vez{{\ve z}}

\begin{document}

\title{Reconstructing biochemical cluster networks}
\author{Utz-Uwe Haus}

\author{Raymond Hemmecke}

\author{Sebastian Pokutta}

\address{ETH Zurich, Switzerland}

\email{utz-uwe.haus@ifor.math.ethz.ch}

\address{Technische Universität Munich, Germany}

\email{hemmecke@ma.tum.de}

\address{University of Erlangen, Germany}

\email{sebastian.pokutta@math.uni-erlangen.de}

\begin{abstract}
Motivated by fundamental problems in chemistry and biology we study
cluster graphs arising from a set of initial states $S\subseteq
\Z^n_+$ and a set of transitions/reactions $M\subseteq\Z^n_+\times\Z^n_+$. The
clusters are formed out of states that can be mutually transformed
into each other by a sequence of reversible transitions. We provide a
solution method from computational commutative algebra that allows for
deciding whether two given states belong to the same cluster as well as for
the reconstruction of the full cluster graph. Using the cluster graph
approach we provide solutions to two fundamental questions: 1) Deciding
whether two states are connected, e.g., if the initial state can
be turned into the final state by a sequence of transition and 2)
listing concisely all reactions processes that can accomplish that. As
a computational example, we apply the framework to the
permanganate/oxalic acid reaction.
\end{abstract}



\maketitle

\section{Introduction}
\label{sec:introduction}
The reconstruction of all biochemical reaction networks composed of a
given finite set of reactions/transitions that explain a given overall
reaction from an initial
state to some final state is one of the fundamental problems in
biology and chemistry alike. The extraction of this decomposition is
essential and has a wide range of applications from the design of
large-scale reactors in process engineering where the presence of
unexpected side products can disrupt the reaction process to the
derivation of rate laws in physical chemistry.

One initial subproblem to be solved is to decide whether there exists
a reaction network at all that explains the given overall reaction. For
example, it has been settled only recently in \cite{CHEMNETS-1}, that
if no additional catalyst is available, the 19 species postulated in
\cite{kovacs2004rmp} do not suffice to explain the
permanganate/oxalic acid reaction. So far, automatic tools to
reconstruct the reaction networks are rare. Even the very promising
recent approach using integer programming \cite{kovacs2004dpo}
suffered from ad-hoc assumptions. Apart from that, there have been
considerable advances in the study of biochemical reaction networks using
methods from discrete mathematics and computer algebra which are not
directly related to our work although they use similar methods (see
e.g., \cite{Faulon2001fp}, \cite{Marwan+Wagler+Weismantel},
\cite{shiu2009scr}).

We will formulate the underlying mathematical problem of computing cluster networks and present a solution approach from commutative algebra. Given two states, represented as vectors $\ve s, \ve t\in \Z^n_+$ (these could be an initial and a final state of a biochemical reaction) and a set of potential transitions $M=U\cup D\subseteq \Z^n_+\times\Z^n_+$, where a transition $\veu\to\vev$ is encoded as the vector $(\veu,\vev)$. The set $U$ represents undirected transitions (that is, with every $(\veu,\vev)\in U$ also $(\vev,\veu)\in U$) whereas $D$ represents directed transitions. A transition $\veu\to\vev$ is applicable at a state $\vea\in\Z^n_+$ if $\vea\geq\veu$. In this case $\vea$ is transitioned to $\vea-\veu+\vev\in\Z^n_+$. We say that $\ve s\in\Z^n_+$ is $M$-connected to $\ve t\in\Z^n_+$ (short: $\ve s \rightarrow_M \ve t$) if there exists a transition path from $\ve s$ to $\ve t$ (in $\Z^n_+$) using only transitions from $M$. We can formulate the following questions:

\bigskip

{\bf Question 1.}
Given two states $\ve s, \ve t\in\Z^n_+$. Decide whether $\ve s \rightarrow_M \ve t$.

\bigskip

{\bf Question 2.}
Given two states $\ve s, \ve t\in\Z^n_+$ such that $\ve s\rightarrow_M \ve t$. Find \emph{all} directed paths from $\ve s$ to $\ve t$ in $\Z^n_+$ using only transitions in $M$.

\bigskip

For small instances, Question 1 and to some extent also Question 2 can be solved using a purely enumerative approach. This approach is finite if we assume that $M$ is homogeneous with respect to some positive grading. This road was successfully pursued in \cite{CHEMNETS-1} to examine the permanganate/oxalic acid reaction. Unfortunately, this approach fails for larger problem sizes due to the vast amount of possibilities that have to be enumerated as a consequence of the {\em combinatorial explosion}. The situation is especially challenging, as due to reversible reactions a large number of reaction paths are similar and thus block the view onto structurally different networks.

In this article, we provide an algorithmic solution to the following two problems. The first problem is the identification of equivalence classes (or \emph{clusters}) of states in $\Z^n_+$ induced by the reversible reactions $U$, that is, given $\ve u, \ve v\in\Z^n_+$, decide whether $\ve u\leftrightarrow_U\ve v$. The second problem is the construction of the graph of all clusters reachable from a given cluster. This cluster graph is the compressed reachability graph starting at $\ve s\in\Z^n_+$ and using only transitions from $M$, where equivalent states with respect to $U$ are contracted.

We will use techniques from computer algebra in order to transform the problems at hand into algebraic ones and then, using Gr\"obner bases, provide two algorithms that solve these problems. We will then propose a solution to Question 1 and Question 2 using the aforementioned cluster approach. The cluster graph that arises from contracting the state graph is usually considerably smaller and both questions can be decided on the cluster graph using traditional graph theoretic methods, thus enabling the successful processing of significantly larger instances. In particular, transitions connecting two clusters are of major importance for any decomposition. They can be easily identified within the cluster graph. Thus, being able to compute cluster graphs, one may study more complex reactions whose symmetries and equivalent paths block the view onto essential parts of the decomposition, see Figure \ref{fig:clustGraph}.

The outline of the article is as follows. As a motivation, we provide an introduction to the network reconstruction problem from a chemical point of view in Section~\ref{sec:motChem}. We will introduce the necessary notation and establish the link to our mathematical approach using directed graphs and commutative algebra. In Section~\ref{sec:preliminaries} the necessary preliminaries and definitions from a mathematical point of view are formulated and basic results are established. We will then derive the described algorithms in Section~\ref{sub:clusterApproach} and provide computational results in Section~\ref{sec:compRes}. Finally, we conclude with some remarks in Section~\ref{sec:conclusion}.

The notation we use is standard (cf.~\cite{cox2007iva,kreuzer2000cca}). These books also provide excellent introductions to commutative computer algebra.

\section{Motivation from chemistry}
\label{sec:motChem}

The decomposition of an overall chemical reaction into elementary reaction steps is one of the fundamental questions in chemistry, since the network of elementary reactions encodes the dynamics of the chemical system and hence allows an analytical examination. We call a reaction {\em elementary} if at most two species react. For example, $2\ce{H2O} \to 2\ce{H2}+\ce{O2}$ is an elementary reaction, but the reverse reaction $2\ce{H2}+\ce{O2}\to 2\ce{H2O}$ is not, as three (not necessarily different) species react.

The decomposition problem can be solved via several steps. We demonstrate them using the well-known permanganate/oxalic acid reaction
\[
  2\ce{MnO4-}+6\ce{H+}+5\ce{H2C2O4} \to 2\ce{Mn^{2+}} + 8\ce{H2O} + 10\ce{CO2}
\]
that has been studied since 1866 \cite{harcourt1866} and that has been studied by several groups (see e.g., \cite{abel1952vrz}, \cite{adler1955mpo}, \cite{Deiss1926},\cite{kovacs2004rmp}, \cite{pimienta1994rmi}, \cite{skrabal1904upo}, \cite{szalkai2000nga}). As a first step, we have to choose a set of possible intermediate species that can be used to explain the overall reaction. For example, the following $19$ species are a commonly accepted set \cite{kovacs2004rmp}:
{\small
\[
  \begin{array}{llll}
    \text{H}_2\text{C}_2\text{O}_4 & \text{H}\text{C}_2\text{O}_4^- & \text{H}^+ & \text{C}_2\text{O}_4^{2-} \\
    \text{Mn}^{2+} & \text{Mn}\text{C}_2\text{O}_4 & \text{Mn}\text{O}_4^- & \text{Mn}\text{O}_2\\
    \text{Mn}^{3+} & \text{C}\text{O}_2 & \text{H}_2\text{O} & [\text{Mn}\text{O}_2,\text{H}_2\text{C}_2\text{O}_4] \\
    \text{C}\text{O}_2^- & [\text{Mn}(\text{C}_2\text{O}_4)]^+ & [\text{Mn}(\text{C}_2\text{O}_4)_2]^- &
   \text{[Mn}\text{C}_2\text{O}_4,\text{Mn}\text{O}_4^-,\text{H}^+\text{]} \\
 \text{[Mn}\text{C}_2\text{O}_4^{2+},\text{Mn}\text{O}_3^-\text{]}^+ &
      \text{[Mn}\text{C}_2\text{O}_4^{2+},\text{Mn}\text{O}_3^-,\text{H}^+\text{]}^{2+} &
   \text{[H}^+,\text{Mn}\text{O}_2,\text{H}_2\text{C}_2\text{O}_4\text{]}&
  \end{array}
\]}

Now, all possible elementary reactions among the postulated species can be computed. For this, note that any chemical state can be represented by a vector $\ves\in\Z^n_+$, where $s_i$ counts how many times species $i$ is present in the state. Moreover, a chemical reaction can be encoded as an integer vector $\ve d=(\ved',\ved'')\in\Z^n_+\times\Z^n_+$. In particular, in component $d'_i$ we specify how many units of species $i$ react and in component $\ved''_i$ we specify how many units of species $i$ are created. In our example below, we will assume for the elementary reactions that no species appears both as reactant and as product, that is, $\supp(\ved')\cap\supp(\ved'')=\emptyset$. Clearly, if we drop this condition we can also model reactions that need some catalyst to be started. Then the catalyst would simply be put on both sides of the reaction in equal amounts.

As any chemical reaction must fulfill a balance of mass and charge, $-\ved'+\ved''$ is a solution of a certain linear system of equations $A\vez=\ve 0$ with $A\in\Z^{m\times n}$ and $\vez\in \Z^n$. For our example, we obtain
\begin{equation}
\label{eq:oxalic}
A = \left(
\begin{array}{rrrrrrrrrrrrrrrrrrr}
0 & 0 &  0 & 0 &   1 & 1 & 1 & 1 &  1 &  0 & 0 & 1 & 0 &  1 &  1 & 2 &  2 &  2 &  1\\
2 & 2 &  0 & 2 &   0 & 2 & 0 & 0 &  0 &  1 & 0 & 2 & 1 &  2 &  4 & 2 &  2 &  2 &  2\\
2 & 1 &  1 & 0 &   0 & 0 & 0 & 0 &  0 &  0 & 2 & 2 & 0 &  0 &  0 & 1 &  0 &  1 &  3\\
4 & 4 &  0 & 4 &   0 & 4 & 4 & 2 &  0 &  2 & 1 & 6 & 2 &  4 &  8 & 8 &  7 &  7 &  6\\
0 & 1 & -1 & 2 &  -2 & 0 & 1 & 0 & -3 &  0 & 0 & 0 & 1 & -1 &  1 & 0 & -1 & -2 & -1\\
\end{array}
\right),
\end{equation}
where the first four rows correspond to the mass balance equations of $\ce{Mn}, \ce{C}, \ce{O}$, and $\ce{H}$, respectively. The last row encodes balance of charge. The $19$ columns correspond to the $19$ postulated species. Thus, by construction, the first $m-1$ rows of this matrix $A$ are always nonnegative and the last row may contain negative entries. The permanganate/oxalic acid reaction is encoded in the (negative and positive parts of the) vector
\[
-\ved'+\ved''=\left(
  \begin{array}{rrrrrrrrrrrrrrrrrrr}
  -5 &  0 & -6 &  0 &  2 &
   0 & -2 &  0 &  0 & 10 &
   8 &  0 &  0 &  0 &  0 &
   0 &  0 &  0 &  0\\
  \end{array}
\right).
\]
The set of elementary reactions that we are looking for is now characterized by all integer vectors $-\ved'+\ved''\in\ker(A)$ with $\supp(\ved')\cap\supp(\ved'')=\emptyset$ and $\|\ved'\|_1\leq 2$. In our example, they can be computed via $\binom{19}{2}+\binom{19}{1}+\binom{19}{1}=209$ linear Diophantine systems. In total, they have $1022$ solutions.

\begin{figure}
\begin{center}
  \includegraphics[height=4cm]{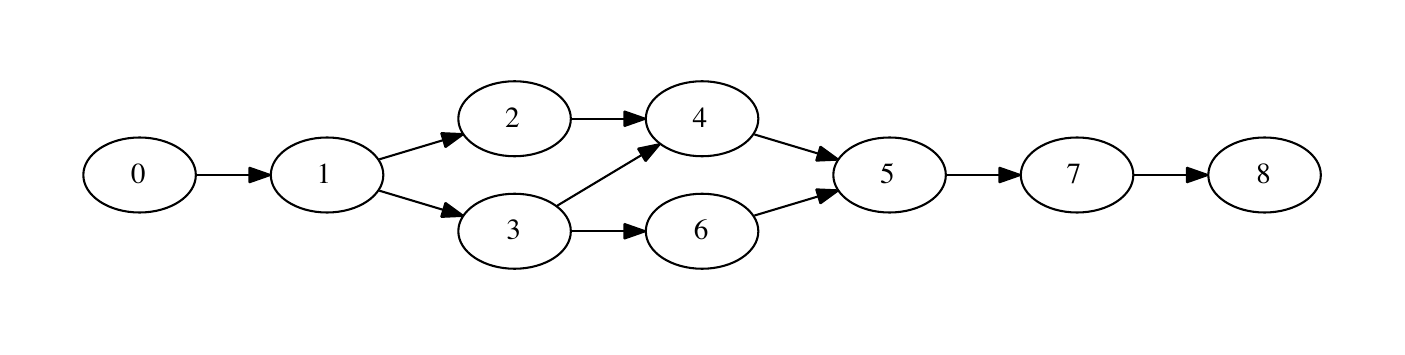}
\end{center}
\caption{Cluster graph for permanganate/oxalic acid reaction. There is a unique reaction connecting cluster $0$ to cluster $1$, showing that it is essential in \emph{any} decomposition.}
\label{fig:clustGraph}
\end{figure}

Finally, these $1022$ reactions can be used to decompose the overall reaction into elementary steps. In particular, we are interested in all such possible reaction networks in order to find a network that explains the observations from experiments most consistently. From a computational point of view, this final decomposition step is by far the most challenging. For example, it was shown in \cite{CHEMNETS-1} that the overall permanganate/oxalic acid reaction cannot be decomposed at all using the $1022$ elementary reactions only. Thus, the $19$ species do not suffice in order to explain this
overall reaction. A solution was given by including a suitable additional species such as $\ce{H2O2}$. Unfortunately, using this additional species and the resulting additional elementary reactions, the reaction network of the overall reaction contains a huge number of paths from the initial state $2\ce{MnO4-}+6\ce{H+}+5\ce{H2C2O4}$ to the final state $2\ce{Mn^{2+}} + 8\ce{H2O} + 10\ce{CO2}$. Many of these paths correspond to identical or essentially identical reaction networks. In order to identify important reactions within the network, we can cluster states together that are connected via a sequence of reversible reactions. Note that we may also cluster according to strongly connected components of states but clustering according to reversible reactions preserves more information of the original state network. In doing so, a large amount of the combinatorial explosion is removed. The resulting coarser network for the example at hand is given in Figure~\ref{fig:clustGraph}.
This cluster network exhibits important elementary reactions needed to explain the overall reaction: For our example it shows that, given the postulated species, at most three essentially different reaction networks exist; some of which might still turn out to be impossible due to
chemical restrictions. Moreover, there is a unique reaction, $\ce{H2C2O4}+\ce{MnO4-}\to\ce{CO2-}+\ce{MnO2}+2\ce{H2O2}+\ce{CO2}$, connecting cluster $0$ to cluster $1$, showing that this reaction must occur in \emph{any} decomposition of the overall chemical reaction.

\section{The cluster graph framework}
\label{sec:preliminaries}

In this section, we want to introduce the necessary notions and definitions. We will also present a few fundamental results that we will use in the following exposition. The objects of our interest will be directed graphs and in particular their connectivity structure.

In the following, we will call the elements in $\Z^n_+$ \emph{states}. Let $U\subseteq\Z^n_+\times\Z^n_+$ be a set of {\em reversible (undirected) transitions} (that is, with every $(\ver',\ver'')\in U$ also $(\ver'',\ver')\in U$), let $D\subseteq\Z^n_+\times\Z^n_+$ be a set of {\em irreversible (directed) transitions} (that is, for every $(\ver',\ver'')\in D$ we have $(\ver'',\ver')\not\in D$), and let $M=U\cup D$ be the set of \emph{transitions}. If $\ves\in\Z^n_+$ is a state and if $(\ver',\ver'')\in M$, then we can transition from state $\ves$ to state $\ves-\ver'+\ver''$ if $\ves\geq\ver'$. Note that if the reversible transition $(\ver',\ver'')\in U$ is applicable at $\ves$, then $(\ver'',\ver')\in U$ is applicable at $\ves-\ver'+\ver''$ leading us back to the original state $\ves$. Using states and transitions, we construct the infinite graph $\Gamma_M$ with node set $\Z^n_+$ and with an arc from $\ve a\in\Z^n_+$ to $\ve b\in\Z^n_+$ if and only if there is some $(\ver',\ver'')\in M$ such that $\ver'\leq\vea$ and such that $\vea-\ver'+\ver''=\veb$ (that is, $\vea$ can be directly transitioned into $\veb$ by some transition from $M$). Given two states $\ves,\vet\in\Z^n_+$, we say that $\ves$ is \emph{$M$-connected} to $\vet$ (short: $\ves\rightarrow_M \vet$) if there exists a directed path in $\Gamma_M$ from $\ves$ to $\vet$.

We assume that $M$ is homogeneous with respect to a positive (multi-) grading $\deg:\Z^n_+\to\Z_+^k$, that is, for all $(\ver',\ver'')\in M$ we have $\deg(\ver')=\deg(\ver'')$. This implies that there are only finitely many states (in $\Z^n_+$) with a given fixed degree. In particular, there are only finitely many states reachable from any given state using the transitions from $M$. We outline this in more detail in the following example.

\begin{example}
Consider the matrix $A\in\Z^{5\times 19}$ given in \eqref{eq:oxalic} which is the mass/charge balance matrix for the permanganate/oxalic acid reaction assuming $19$ species.

Adding up the first $4$ rows we obtain the vector:
\[
  a_{1}+a_{2}+a_{3}+a_{4} =
  \left(
    \begin{array}{rrrrrrrrrrrrrrrrrrr}
      8 &  7 &  1 &  6 & 1 &  7 &  5 &  3 & 1 &  3 &  3 & 11 &  3 &  7 & 13 & 13 & 11 & 12 & 12\\
    \end{array}
  \right),
\]
which gives a desired positive grading. We can in fact also obtain a positive multi-grading by adding the vector $a_{1}+a_{2}+a_{3}+a_{4}$ sufficiently often to the rows of $A$:
\[
  \bar{A} =
  \left(
    \begin{array}{rrrrrrrrrrrrrrrrrrr}
     8 &  7 &  1 &  6 & 1 &  8 &  6 &  4 & 2 &  3 &  3 & 12 &  3 &  8 & 14 & 15 & 13 & 14 & 13\\
    10 &  9 &  1 &  8 & 1 &  9 &  5 &  3 & 1 &  4 &  3 & 13 &  4 &  9 & 17 & 15 & 11 & 14 & 14\\
    10 &  8 &  2 &  6 & 1 &  7 &  5 &  3 & 1 &  3 &  5 & 13 &  3 &  7 & 13 & 14 & 11 & 13 & 15\\
    12 & 11 &  1 & 10 & 1 & 11 &  9 &  5 & 1 &  5 &  4 & 17 &  5 & 11 & 21 & 21 & 18 & 19 & 18\\
    32 & 29 &  3 & 22 & 2 & 28 & 19 & 12 & 1 & 12 & 12 & 44 & 13 & 27 & 53 & 52 & 43 & 46 & 47\\
    \end{array}
  \right)
\]
Then the multi-grading is obtained by setting $\deg(\vea):=\bar{A}\vea$ for $\vea\in\Z^n_+$.
\end{example}

Now let us partition the graph $\Gamma_M$ into clusters. This approach
was initially suggested in \cite{CHEMNETS-1}. However, in contrast
to partitioning via the strongly connected components of $\Gamma_M$ as
in \cite{CHEMNETS-1}, we partition via equivalence classes which group
all states $\ve v$
that can be reached from a specific $\ve u$ using only transitions
from $U$. These equivalence classes will be called
\emph{clusters}. Note, that $\ve u \rightarrow_U \ve v$ if and only if
$\ve v \rightarrow_U \ve u$. Thus the relation $\ve u \rightarrow_U
\ve v$ is reflexive, symmetric, and transitive and therefore indeed an
equivalence relation on $\Gamma_M$ and we write $\ve u
\leftrightarrow_U \ve v$. By $\cluster(\ve u)\subseteq\Z^n_+$ we
denote the cluster of $\ve u \in \Z^n_+$ which is defined in the
canonical way, i.e.,
$$\cluster(\ve u) := \{\ve v \in \Z^n_+ \mid \ve u \leftrightarrow_U
\ve v\}.$$

Now, the remaining set $D\subseteq\Z^n$ of irreversible transitions
defines a directed graph with the clusters (of $\Z^n_+$) as vertices and
with a directed edge from cluster $\cluster(\ve u)$ to cluster
$\cluster(\ve v)$ if there exist two states $\ve x \in \cluster(\ve u)$
and $\ve y \in \cluster(\ve v)$ such that $\vex$ can be transitioned into $\vey$ via some $\ved\in D$. We denote
this cluster graph by $\clustergraph(U,D)$. Be aware that we indeed
only add the arc $(\cluster(\ve u), \cluster(\ve v))$ if the cluster
$\cluster(\ve v)$ can be reached from the cluster $\cluster(\ve u)$ by
a transition path of length $1$.

In applications, we are usually given a finite set $S$ of initial
states. Note that in our chemical setting, $S$ typically contains only
one state. More than one initial state could occur if a set of
experiments (i.e., pairs of initial and final states) have to be explained
(by networks). The nodes of the subgraph $\bar \Gamma = (\bar
S, \bar M)$ of $\Gamma_M$ reachable from $S$ in $\Gamma_M$ are
\[
  \bar S = \left\{ \ves + \sum_{i = 1}^{k} (-\veu_i+\vev_i)  \mid k\in \mathbb
    N, \ve s \in S, \{(\veu_i,\ve v_i)\mid i\in [k]\} \subseteq M \text{
      s.t. } \ve s + \sum_{i = 1}^{l-1} (-\veu_i+\vev_i)\geq \veu_l\ \forall l \in
    [k]\right\}.
\]
Herein, we use for convenience $[k]:=\{1, \dots, k\}$ for $k\in\Z_+$. As the arc set $\bar M\subseteq\bar S\times\bar S$ we obtain the set of those arcs in $\Gamma_M$ involving only nodes from $\bar S$. The important difference to a classical directed graph is now that $\bar M$ and $\bar S$ are not given explicitly but implicitly by a set of transitions/reactions $M\subseteq\Z^n_+\times\Z^n_+$ and a set of initial states $S\subseteq\Z^n_+$. We call $\bar \Gamma = (\bar S, \bar M)$ the {\em state graph} of $S$ and with respect to $M$. As $\bar S$ and $\bar M$ might be already very large for small instances it is favorable to not completely calculate $\bar S$ and $\bar M$. In these cases, computing the much smaller cluster graph may still give important information on the decomposition of the overall reaction.

As $M$ is assumed to be homogeneous with respect to some (positive) grading, the state graph of $S$ w.r.t. M decomposes for $\ves_1, \ves_2$ with different degrees and thus one can confine the analysis to sets $S$ with elements of the same degree.

Having provided the considered setting, we set out to provide algorithmic solutions to answer Questions 1 and 2 using cluster graphs. For this we show how to construct the part $\clustergraph(U,D,S)$ of the cluster graph $\clustergraph(U,D)$ reachable from some $\cluster(\ve s)$, $\ve s\in S$.

\section{Reconstructing cluster graphs}
\label{sub:clusterApproach}
In this section we will use computer algebraic tools to reconstruct the cluster graph reachable from states given in a finite set $S$ via transitions in $M=U\cup D$. The positive grading implies that each cluster contains only finitely many states. Moreover, although the total cluster graph over $\Z^n_+$ has infinitely many clusters as vertices, the positive grading implies that for any given state $\ves\in\Z^n_+$ only finitely many clusters can be reached from the clusters $\cluster(\ves)$ with $\ve s\in S$, within the cluster graph. In order to reconstruct $\clustergraph(U,D,S)$, we have to solve the following two problems:

\bigskip

{\bf Main Problem 1.}
Given two states $\ves, \vet\in\Z^n_+$. Decide whether $\cluster(\ves)=\cluster(\vet)$, i.e., whether $\ves\leftrightarrow_U \vet$.

\bigskip

{\bf Main Problem 2.}
Given $\ve s\in\Z^n_+$. List all transitions $\ved=(\ved',\ved'')\in D$ that are applicable to at least one state $\ve s_{\ve d}\in\cluster(\ve s)$,
that is, $\ves_{\ve d}\geq\ved'$.

\bigskip

Main Problem 1 captures the problem of being able to decide whether
two states are in the same equivalence class whereas Main Problem 2
has to be solved in order to identify clusters reachable from a given
cluster. We will provide a solution for both problems in form of
computer algebraic algorithms. If the clusters are small enough such
that all states can be enumerated explicitly, Main Problem 1 can be solved
by simple enumeration. However, the sizes of the clusters may prohibit
an explicit enumeration. Therefore an approach is sought that does not
explicitly enumerate the cluster elements. We will tackle this problem
by transforming the set $U$ into a binomial ideal. We consider the
polynomial ring $K[X]$ where $K$ is an arbitrary field and $X$ is the
set of the variables. Let $J_U$ be defined to be the ideal
\[
  J_U:=\left\langle \ve x^{\veu}-\ve x^{\vev}: (\veu,\vev)\in U\right\rangle \subseteq K[X].
\]
Then $\ve y \leftrightarrow_U \ve z$ if and only if $x^{\ve y} -
x^{\ve z} \in J_U$ as shown in the following theorem. This theorem has been rediscovered many times --- an early reference is \cite{mayr1982cwp}. For the sake of completeness we include a proof below (cf.\ \cite{Hemmecke2009}).
\begin{theorem}\label{thm:testConnect}
  Let $U\subseteq\Z^n_+\times\Z^n_+$ and $\vey,\vez\in\Z^n_+$. Then $\vey\leftrightarrow_U \vez$ if and only if $\vex^{\vey}-\vex^{\vez}\in J_U$.
\begin{proof}
  Let $\vey,\vez\in\Z^n_+$ such that $\ve y \leftrightarrow_U \vez$. Thus, there exists a path from $\ve y$ to $\ve z$ in $\Z^n_+$. More specifically, there exists a sequence of points $(\vep_1,\vep_2,...,\vep_k)\subseteq\Z^n_+$ where $\vep_1=\vey$, $\vep_k=\vez$, and $\vep_i-\vep_{i+1}=\veu'_i-\veu''_i$ for some $(\veu'_i,\veu''_i)\in U$ for $i=1,...,k-1$. Thus, there exists $\ve\gamma_i\in\Z^n_+$ such that
  $\ve p_i=\veu'_i+\ve\gamma_i$, and $\ve p_{i+1}=\veu''_i+\ve\gamma_i$ for every $i=1,...,k-1$. Hence, $\vex^{\vep_i}-\vex^{\vep_{i+1}}=\vex^{\ve\gamma_i}(\vex^{\veu'_i}-\vex^{\veu''_i})$, and therefore,
  \[
    \vex^{\vey}-\vex^{\vez}=\sum_{i=1}^{k-1} (\vex^{\vep_i}-\vex^{\vep_{i+1}})=\sum_{i=1}^{k-1}\vex^{\ve\gamma_i}(\vex^{\veu'_i}-\vex^{\veu''_i})\in J_U
  \]
  as required.

  Conversely, assume that $\ve x^{\ve y}-\ve x^{\ve z}\in J_U$. Further, suppose for contradiction, $\ve y \;\not\hspace{-0.5ex} \leftrightarrow_U \ve z$.
  As $\vex^{\vey}-\vex^{\vez}\in J_U$, we may write $\vex^{\ve y}-\vex^{\ve z}=\sum_{i=1}^dc_i\vex^{\ve\gamma_i}(\vex^{\veu'_i}-\vex^{\veu''_i})$ where $(\veu'_i,\veu''_i)\in U$, $c_i\in K$, and $\ve\gamma_i\in\Z^n_+$. Note that we allow $(\veu'_i,\veu''_i)=(\veu'_j,\veu''_j)$ for $i\ne j$. Now, let $I := \{i \in [d] \mid (\ve\gamma_i+\veu''_i) \leftrightarrow_U \vey\}$. Clearly, $(\ve\gamma_i+\veu''_i) \in \Z^n_+$ for all $i\in [d]$. Note that if $(\ve\gamma_i+\veu''_i) \leftrightarrow_U \vey$ then $(\ve\gamma_i+\veu'_i) \leftrightarrow_U \vey$ since $\ve\gamma_i+\veu''_i\geq\veu''_i$ and $(\ve\gamma_i+\veu''_i)+(-\veu''_i+\veu'_i)=(\ve\gamma_i+\veu_i')$. Thus, the set of monomials consisting of $\vex^{\ve\gamma_i}\vex^{\veu''_i}$ and $\vex^{\ve\gamma_i}\vex^{\veu'_i}$ for all $i\in I$, which includes $\vex^{\ve y}$ and not $\vex^{\vez}$, is disjoint from the set of monomials consisting of $\vex^{\ve\gamma_i}\vex^{\veu''_i}$ and $\vex^{\ve\gamma_i}\vex^{\veu'_i}$ for all $i\not\in I$, which includes $\vex^{\vez}$ and
  not $\vex^{\vey}$. Let $f(\vex) = \sum_{i\in I}c_i\vex^{\ve\gamma_i}(\vex^{\veu''_i}-\vex^{\veu'_i})$ and let $g(\vex) = -\sum_{i\not\in I}c_i\vex^{\ve\gamma_i}(\vex^{\veu''_i}-\vex^{\veu'_i})$. It is readily seen that the polynomials $f(\vex)$ and $g(\vex)$ have a disjoint set of
  monomials, and therefore, $f(\vex)=\vex^{\ve y}$ and $g(\vex)=\vex^{\vez}$ since $\vex^{\vey}-\vex^{\vez}=f(\vex)-g(\vex)$. However, this is impossible since $f(\1)=0$ and $g(\1)=0$ but $\1^{\vey}=1$ and $\1^{\vez}=1$.
\end{proof}
\end{theorem}

Using Theorem~\ref{thm:testConnect} we can now solve Main Problem 1 as follows. Choose $\prec$ to be an arbitrary term ordering, let $G_\prec(J_U)$ be a Gr\"obner basis of $J_U$ with respect to $\prec$, and let $\ve y,\ve z\in \Z^n_+$ be two states. Then $\vey\leftrightarrow_U \ve z$ if and only if $\ve x^{\ve y}-\ve x^{\vez}\in J_U$ by Theorem~\ref{thm:testConnect}. As $G_\prec(J_U)$ is a Gr\"obner basis, $\ve x^{\ve y}-\ve x^{\ve z}\in J_U$ if and only if $\NF_\prec(\ve x^{\ve y}-\ve x^{\ve z},G)=0$, where $\NF_\prec(.)$ is the normal form operator. Note that once the Gr\"obner basis $G_\prec(J_U)$ of $J_U$ has been calculated, the membership test $\NF_\prec(\ve x^{\ve y}-\ve x^{\vez},G)=0$ can be performed easily.

As an immediate consequence of this construction, it follows that
every cluster $\cluster(\ve s)$ with $\ve s\in \Z^n_+$ has a unique
($\prec$-minimal) representative outside the leading term ideal
$\lt_\prec(J_U)$ which can be obtained by calculating the normal
form. Slightly abusing notation we denote this unique minimal element
by $\cluster(\ve s) := \NF_\prec(\ve x^{\ve s},G)$. Thus, the
monomials outside of $\lt_\prec(J_U)$ are in one-to-one correspondence
with possible clusters in $\Z^n_+$. Consequently, the number of
clusters reachable from $\cluster(\ve s)$ is bounded by the number of
monomials outside of $\lt_\prec(J_U)$ which have the same (multi-)
degree $\deg(\ve s)$. These numbers, however, are encoded in the
multi-graded Hilbert-Poincar\'e series \cite{kreuzer2000cca} of $J_U$:
\[
  \HPS_{(\vez)}:=\HPS(\vez,J_U,\prec,\deg):=\sum_{\ve\alpha\in\Z^n_+:x^{\ve\alpha}\not\in\lt_\prec(J_U)}\vez^{\deg(\ve\alpha)}.
\]
Note that this potentially infinite sum can always be written as a
rational function $p(\ve z)/q(\ve z)$. Moreover, it does not depend on
the actual term ordering $\prec$ chosen and it can be computed
from any Gr\"obner basis of $J_U$ rather efficiently. The
coefficient in front of $\ve z^{\deg(\ve s)}$ in the multi-variate
Hilbert series expansion of $p(\ve z)/q(\ve z)$ is exactly the number
of clusters of degree $\deg(\ve s)\in\Z^k$. Thus, we can compute in
advance an upper bound on the number of clusters that can be reached
from $\cluster(\ve s)$ or that can reach $\cluster(\ve s)$ inside the
cluster graph $\clustergraph(U,D)$.

Recall that for a directed transition $(\ved',\ved'')\in D$ and a state $\vez\in\Z^n_+$ the transition $(\ved',\ved'')$ is applicable at $\vez$, if $\vez\geq\ved'$. Thus, in order to solve Main Problem 2, we need to decide whether a given directed transition $(\ved',\ved'')\in\Z^n_+\times\Z^n_+$ is applicable at some state $\vez\in\cluster(\ve y)$ given only some representative $\ve y$ of the cluster $\cluster(\ve y)$. The following theorem will be a main ingredient in order to solve this question. Recall that for an ideal $I\subseteq K[X]$ and a vector $\vem\in\Z^n_+$ the colon ideal $I:{\vex}^{\ve m}$ is defined as $I:{\ve x}^{\ve m} :=\{p\in K[X] \mid p {\vex}^{\vem}\in I\}$.
\begin{theorem}\label{Theorem: Connection above d}
  Let $U\subseteq\Z^n$, $\bar{\ve d}\in\Z^n_+$, and $\ve y, \ve z\in\Z^n_+$ with $\ve y, \ve z\geq\bar{\ve d}$. Moreover, let $V\subseteq\Z^n$ such that $J_V=J_U:\ve x^{\bar{\ve d}}$. Then
  \[
    \ve y\leftrightarrow_U\ve z\;\;\; \Leftrightarrow\;\;\; \ve y-\bar{\ve d}\leftrightarrow_V\ve z-\bar{\ve d}.
  \]
\begin{proof}
  This follows immediately using Theorem \ref{thm:testConnect}: Note that $\vey\leftrightarrow_U\vez$ if and only if $\vex^{\vey}-\vex^{\ve z}\in J_U$ and similarly, $\ve y-\bar{\ved}\leftrightarrow_V\ve z-\bar{\ve d}$ if and only if $\ve x^{\vey-\bar{\ve d}}-\ve x^{\ve z-\bar{\ve d}}\in J_V$. Now observe that $\vex^{\vey-\bar{\ved}}-\vex^{\vez-\bar{\ved}}\in J_V$ holds if and only if $\vex^{\vey}-\vex^{\vez}\in J_U$, as $J_V=J_U:\vex^{\bar{\ved}}=\left\langle\vex^{\veu}-\vex^{\vev}:\veu+\bar{\ve d}\leftrightarrow_U \vev+\bar{\ved}\right\rangle$.

  Observe that for the backward implication of the equivalence, we employ the condition $\vey, \vez\geq\bar{\ved}$, that is, $\vex^{\vey-\bar{\ved}}-\vex^{\vez-\bar{\ved}}$ is indeed a polynomial in $\ve x$.
\end{proof}
\end{theorem}

\begin{figure}
  \begin{center}
   \scalebox{0.9}{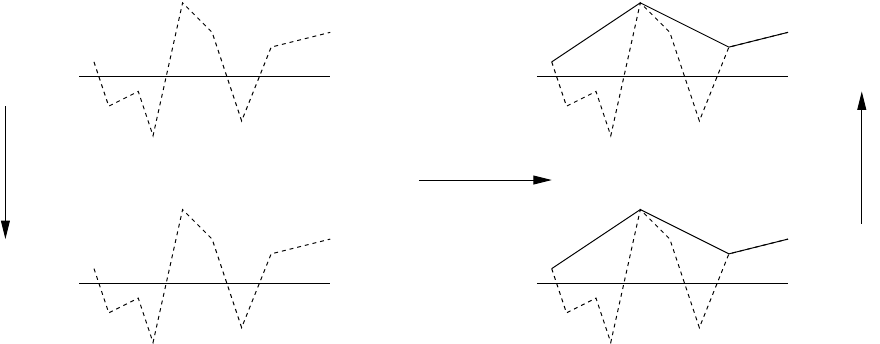}
  \end{center}
  \caption{Shift relative to $\bar{\ve d}$.}
  \label{fig:saturation}
\end{figure}

Note that Theorem \ref{Theorem: Connection above d} shows that $V$ connects any two points $\vey, \vez$ in $\cluster(\vey)$ with $\vey,\vez\geq\bar{\ve d}$ via a sequence of points in $\cluster(\ve y)$ whose coordinates are also at least as big as $\bar{\ve d}$. Such a path can be constructed by shifting the relation $\ve y-\bar{\ved}\leftrightarrow_V\ve z-\bar{\ve d}$ by $\bar{\ved}$ as shown in Figure~\ref{fig:saturation}. Moreover, by computing the normal form of $\ve y$ with respect to a $\prec$-Gr\"obner basis of $J_V$ where $\prec$ is a term ordering that maximizes
the $j$-th component, we can decide for a given $\ve d:=\bar{\ve d}+d_j\ve e_j$ whether there is some $\ve z\in\cluster(\ve y)$ with $\vez\geq\ve d$. This is summarized in Algorithm \ref{alg:CI}.

\restylealgo{boxed}\linesnumbered
\incmargin{1em}
\begin{algorithm}
\label{alg:CI}
\caption{Coordinate increment algorithm (CI).}
\SetKwFunction{Compute}{compute}
\SetKwFunction{Set}{set}
\SetKwData{AFalse}{False}
\SetKwInOut{Input}{Input}
\SetKwInOut{Output}{Output}
\Input{Let $U\subseteq\Z^n$, $\ve y, \ve d\in\Z^n_+$. Moreover, let
  $S=\supp(\ve d)$, $j\in S$, and $\bar{\ve d}:=\ve d-d_j\ve e_j$.}
\Output{$\ve z$ if $\ve z\in\cluster(\ve y)$ with $z_j\geq d_j$ exists, \AFalse otherwise.}
\BlankLine
\Set $\bar{\ve  d} := {\ve  d}- d_je_j$;\\
\Compute $V\in\Z^n$ such that
$      J_V=J_U: \ve x^{\bar{\ve d}}
         =\left\langle \ve x^{\ve u}-\ve x^{\ve v}:
             \ve u+\bar{\ve d}\leftrightarrow_U \ve v+\bar{\ve d} \right\rangle;$\\
\Compute a Gr\"obner basis $G$ of $J_V$ w.r.t.~a term ordering $\prec$ that maximizes the $j$-th
    component; \\
$\ve x^{\bar{\ve z}} := \NF_\prec(\ve x^{\ve y-\bar{\ve d}}, G);$ \\
\If{$\bar{z}_j\geq d_j$}{
	{
	$\ve z:=\bar{\ve z}+\bar{\ve d}\in\cluster(\ve y)$; \\
	\Return $\ve z$;
	}\\
	{\Else{\Return \AFalse;}}
}
\end{algorithm}

We will now formulate Algorithm~\ref{alg:cct} that solves Main Problem 2 by iteratively calling the coordinate increment algorithm (Algorithm
\ref{alg:CI}). Algorithm~\ref{alg:cct} constructs an element $\vez\in\cluster(\vey)$ with $\ve z\geq\ved$ for a given $\vey\in\mathbb Z^n_+$ and $\ved\in\Z^n_+$ if such an element $\ve z$ exists.

\restylealgo{boxed}\linesnumbered
\incmargin{1em}
\begin{algorithm}
\label{alg:cct}
\caption{Cluster connectivity test (CCT).}
\SetKwFunction{Compute}{Compute}
\SetKwData{AFalse}{False}
\SetKwFunction{Set}{set}
\SetKwFunction{CI}{CI}
\SetKwFunction{Stop}{STOP}
\SetKwFunction{Let}{let}
\SetKwInOut{Input}{Input}
\SetKwInOut{Output}{Output}
\Input{Let $U, \ve y, d$ be as in Theorem~\ref{Thm:cct}.}
\Output{$\ve z\in\cluster(\ve y)$ with $\vez\geq\ved$ if exists, \AFalse otherwise.}
\BlankLine
\Set $S:= \supp(\ved)$; \\
\Let $\{s_1, \dots, s_k\} := S$ be an enumeration of $S$ with $k\leq n$; \\
\Set $\ved^{(1)} := d_{s_1}\ve e_{s_1}$; \\
\Set $\vez^{(0)} := \vey$; \\
\For{$i \in [k]$}{
	\Set $\ve z^{(i)} := \CI(U,\ve z^{(i-1)},\ve d^{(i)},s_i)$; (with CI from Theorem~\ref{Thm:cct})\\
	\If{$\ve z^{(i)} == \AFalse$}{\Return \AFalse;\\ \Stop;}
	\If{$i<k$}{\Set ${\ve d}^{(i+1)} := {\ve d}^{(i)}+ d_{s_{i+1}}\ve e_{s_{i+1}}$;}
}
\Return $\ve z^{(k)}$;
\end{algorithm}

\begin{theorem}\label{Thm:cct}
  Let $U\subseteq\Z^n$, $\ve y \in \Z^n_+$, and $\ved\in\Z^n$. Then Algorithm~\ref{alg:cct} decides whether there exists $\ve z\in\cluster(\ve y)$ with $\vez\geq\ved$ and if the answer is affirmative it returns $\ve z$.
\begin{proof}
  Let $S$, $\ved^{(i)}$, $\vez^{(i)}$, and $k$ be as calculated by Algorithm~\ref{alg:cct}. Suppose that $\vez\in\cluster(\vey)$ with $\vez\geq\ved$ exists. We show that at the end of iteration $i$ we obtain an element $\vez^{(i)}\in\cluster(\vey)$ with $z^{(i)}_{s_j}\geq d_{s_j}$ for all $j\in [i]$. Then for $i=k=|\supp(\ved)|$ the assertion follows as $z^{(k)}_{s_j}\geq d_{s_j}$ for all $j\in [k]$ and hence $\vez^{(k)}\geq\ved$.

  We have to verify that in iteration $i$ we can apply Algorithm~\ref{alg:CI} to the input $U$, $\vez^{(i-1)}$, $\ved^{(i)}$, and $s_i$. At the start of iteration $i$ we have $\vez^{(i-1)}\geq{\ved}^{(i)}-d_{s_{i}}\vece_{s_{i}}$. In the case of $i=1$ this is clear as ${\ved}^{(1)}-d_{s_{1}}\vece_{s_{1}}=\ve 0$ and $\ve z^{(0)} = \ve y\geq\ve 0$. For $i>1$ this follows by induction as $\ve z^{(i-1)} := \text{CI}(U,\ve z^{(i-2)},\ve d^{(i-1)},s_{i-1})$ with $\ve z^{(i-1)} \geq \ve d^{(i-1)}$ if such an element $\ve z^{(i-1)} \in \cluster(\ve y)$ exists. With $\ve d^{(i-1)} = \ve d^{(i)} - d_{s_{i}}\ve e_{s_{i}}$ the claim follows. Thus we can indeed apply Algorithm~\ref{alg:CI} to the input $U$, $\ve z^{(i-1)}$, $\ve d^{(i)}$, $s_i$ and it returns $\ve z^{(i)}$ with $\ve z^{(i)} \geq \ve d^{(i)}$ if such an element $\ve z^{(i)} \in \cluster(\ve y)$ exists, so that the loop is well defined and we obtain $\ve z^{(i)} \geq \ve d^{(i)}$ at the end of each iteration. Let $i \in [k]$ and observe that $\ve d^{(i)}_{s_j} \geq d_{s_j}$ for all $j\in [i]$ so that, together with $\ve z^{(i)} \geq \ve d^{(i)}$, it follows $z^{(i)}_{s_j}\geq d_{s_j}$ for all $j\in [i]$.

  In case there is no such $\ve z \in\cluster(\ve y)$ with $\vez\geq\ved$, there exists $i\in [k]$ such that $\text{CI}(U,\ve z^{(i-1)},\ve d^{(i)},s_i)$ returns `False' and the algorithm stops. This finishes the proof.
\end{proof}
\end{theorem}

\begin{remark}
  It should be noted that in our chemical setting we have $||\ved'||_1\leq 2$ for any elementary reaction $(\ved',\ved'')$. Therefore, we only have to compute generators for $J_U:(x_ix_j)=(J_U:x_i):x_j$. This square number of Gr\"obner bases can in fact be pre-computed before we start construction of the cluster graph $\clustergraph(U,D,S)$.

  Computation of $I:\vex_i$ for a positively graded binomial ideal $I$ can be done without any additional indeterminate. (Thus, the computation is much faster.) It suffices to compute any DegRevLex-Gr\"obner basis of $I$ with an ordering such that $x_i$ is minimal, and then dividing any binomial by $x_i$ for which this is possible \cite{Hosten+Sturmfels:GRIN}. Note that the resulting set of binomials is still a DegRevLex-Gr\"obner basis of $I:\vex_i$.

  Finally, in order to compute $J_U:(x_ix_j)=(J_U:x_i):x_j$, we first compute a generating set of $J_U:x_i$ as explained in the previous paragraph. As this generating set is a Gr\"obner basis of $J_U:x_i$, we might even use a Gr\"obner walk \cite{Collart+Kalkbrener+Mall} from there to compute the desired DegRevLex-Gr\"obner basis of $J_U:x_i$ with $x_j$ smallest.
\end{remark}

\begin{example}
  With the help of Algorithm~\ref{alg:cct} we may decide for any given (representative state) $\vey\in\Z^n_+$ whether a given directed transition $(\ved',\ved'')\in D$ is applicable at some state $\ve z\in\cluster(\ve y)$. Consequently, cluster $\cluster(\ve y)$ is connected to cluster $\cluster(\vez-\ved'+\ved'')$ in the cluster graph $\clustergraph(U,D)$. However, as the following simple example shows, there may be more than one cluster in $\clustergraph(U,D)$ that can be reached from $\cluster(\ve y)$ via the transition $(\ved',\ved'')\in D$. The simple reason for this is that the state $\ve z\in\cluster(\ve y)$ with $\vez\geq\ved'$ is not unique and a different state $\vez'\in\cluster(\vey)$ with $\vez'\geq\ved'$ could be transitioned to a state $\vez'-\ved'+\ved''$ lying in a different cluster than $\vez-\ved'+\ved''$.

\begin{center}
\scalebox{0.5}{  \begin{picture}(0,0)%
\includegraphics{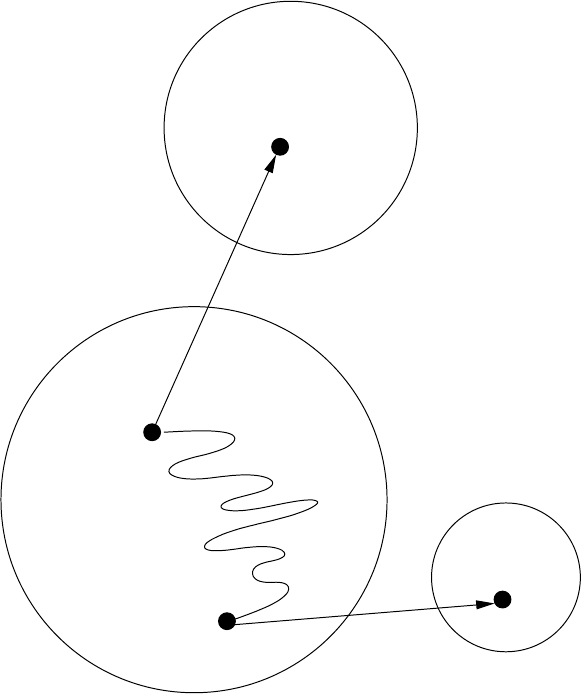}%
\end{picture}%
\setlength{\unitlength}{4144sp}%
\begingroup\makeatletter\ifx\SetFigFont\undefined%
\gdef\SetFigFont#1#2#3#4#5{%
  \reset@font\fontsize{#1}{#2pt}%
  \fontfamily{#3}\fontseries{#4}\fontshape{#5}%
  \selectfont}%
\fi\endgroup%
\begin{picture}(4431,5285)(1226,-6329)
\put(4111,-5926){\makebox(0,0)[lb]{\smash{{\SetFigFont{12}{14.4}{\rmdefault}{\mddefault}{\updefault}{\color[rgb]{0,0,0}$(\ve d',\ve d'')$}%
}}}}
\put(2686,-3166){\makebox(0,0)[lb]{\smash{{\SetFigFont{12}{14.4}{\rmdefault}{\mddefault}{\updefault}{\color[rgb]{0,0,0}$(\ve d',\ve d'')$}%
}}}}
\put(3121,-1996){\makebox(0,0)[lb]{\smash{{\SetFigFont{12}{14.4}{\rmdefault}{\mddefault}{\updefault}{\color[rgb]{0,0,0}$\vex^{\ve z - \ve d' + \ve d''}$}%
}}}}
\put(4726,-5401){\makebox(0,0)[lb]{\smash{{\SetFigFont{12}{14.4}{\rmdefault}{\mddefault}{\updefault}{\color[rgb]{0,0,0}$\vex^{\ve z' - \ve d' + \ve d''}$}%
}}}}
\put(1951,-4576){\makebox(0,0)[lb]{\smash{{\SetFigFont{12}{14.4}{\rmdefault}{\mddefault}{\updefault}{\color[rgb]{0,0,0}$\vex^{\ve z}$}%
}}}}
\put(2566,-6091){\makebox(0,0)[lb]{\smash{{\SetFigFont{12}{14.4}{\rmdefault}{\mddefault}{\updefault}{\color[rgb]{0,0,0}$\vex^{\ve z'}$}%
}}}}
\end{picture}%
}
\end{center}

  Let $\vez=(1,0,1),\vez'=(0,1,1),\veu'=(1,0,1),\veu''=(0,1,1),\ved'=(0,0,1),\ved''=(1,0,0)$, and set $U=\{(\veu',\veu'')\}$ and $D=\{(\ved',\ved'')\}$. By construction, $\vez\leftrightarrow_U\vez'$ and $\vez,\vez'\geq\ved'$, that is, the transition $(\ved',\ved'')\in D$ can be applied to both $\vez$ and $\vez'$. The transitioned states are $\vez-\ved'+\ved''=(2,0,0)$ and $\vez-\ved'+\ved''=(1,1,0)$. These two states do not belong to the same cluster with respect to $U$, since $x_1^2-x_1x_2\not\in J_U=\langle x_1x_3-x_2x_3\rangle$, since every polynomial in $J_U$ is a multiple of $x_3$.
\end{example}

\begin{lemma}
\label{Lemma: Unique neighboring cluster}
  Let $\vey\in\Z^n_+$ be a state, let $(\ved',\ved'')\in D\subseteq\Z^n_+\times\Z^n_+$ and let $U\subseteq\Z^n_+\times\Z^n_+$ be the set of reversible transitions. Moreover, assume that $J_U:\vex^{\ved'}=J_U$. Then for any $\vez,\vez'\in\cluster(\vey)$ with $\vez,\vez'\geq\ved'$, the (transitioned) states $\vez-\ved'+\ved''$ and $\vez'-\ved'+\ved''$ belong to the same cluster with respect to $U$.
\begin{proof}
  As $\vez\leftrightarrow_U\vez'$, we get $\vex^{\vez}-\vex^{\vez'}\in J_U$ by Theorem \ref{thm:testConnect}. Moreover, as $\vez-\ved'\geq\ve 0$ and $\vez'-\ved'\geq\ve 0$, we have that $\vex^{\vez-\ved'}-\vex^{\vez'-\ved'}\in J_U:\vex^{\ved'}$, by definition of $J_U:\vex^{\ved'}$. Thus, we also have $\vex^{\vez-\ved'+\ved''}-\vex^{\vez'-\ved'+\ved''}=\vex^{\ved''}\left(\vex^{\vez-\ved'}-\vex^{\vez'-\ved'}\right)\in J_U:\vex^{\ved'}$. Since $J_U:\vex^{\ved'}=J_U$ by assumption, we conclude by Theorem \ref{thm:testConnect} that the states $\vez-\ved'+\ved''$ and $\vez'-\ved'+\ved''$ belong to the same cluster with respect to $U$.
\end{proof}
\end{lemma}

\begin{remark}
  Note that Lemma \ref{Lemma: Unique neighboring cluster} implies that for every given state $\vey\in\Z^n_+$, there is at most one cluster in $\clustergraph(U,D)$ that is reachable from $\cluster(\vey)$ via any given $(\ved',\ved'')\in D$.
\end{remark}

Now we are finally ready to state an algorithm that reconstructs the cluster graph $\clustergraph(U,D,S)$ when $J_U$ and $D$ satisfy certain conditions.

\begin{lemma}
\label{lem:cgr}
  Let $\ve s \in \Z^n_+$ be an initial state and $M=U\cup D\subseteq \Z^n_+\times\Z^n_+$ be as above. Moreover, assume that $J_U:\vex^{\ved'}=J_U$ for all $(\ved',\ved'')\in D$. Then Algorithm~\ref{alg:reconstruction} reconstructs the induced cluster graph $\clustergraph(U,D,S)$.
\begin{proof}
  The list $N$ contains all the nodes that have been discovered but not visited yet. The list $V$ contains all visited nodes and $E$ is the edge set that we construct successively. Clearly Algorithm~\ref{alg:reconstruction} is finite if the cluster graph $\clustergraph(U,D,S)$ is finite: In each iteration of the loop starting in $5$ we remove one element from $N$ and we only add nodes to $N$ (lines 15/16) if the constructed node is new. It remains to show that the algorithm reconstructs the cluster graph $\clustergraph(U,D,S)$. We initialize $N$ with $\NF(\ve s,G_U)$ and we successively process all the nodes in $N$. For each node $\veu \in N$ we have to decide whether a directed transition $(\ved',\ved'')$ can lead to a new adjacent cluster. (By Lemma \ref{Lemma: Unique neighboring cluster}, there exists at most one such cluster due to our assumption $J_U:\vex^{\ved'}=J_U$.) Therefore we employ Algorithm~\ref{alg:cct} to construct $\vev:= \text{CCT}(U,\ve u,\ved')$ if such a $\vev$ exists. In a next step we verify that the transitioned state $\vev-\ved'+\ved''\not\in\cluster(\ve u)$. In this case we found a transition $(\ved',\ved'')$ leading to an adjacent cluster and we add the corresponding edge $(\ve u, \ve w)$ to $E$. Finally we test if the canonical cluster representative $\ve w:=\NF(\vev-\ved'+\ved'',G_U) \not \in V$, i.e., we discovered a new cluster and add $\ve w$ to $V$ in this case.
\end{proof}
\end{lemma}

\restylealgo{boxed}\linesnumbered
\incmargin{1em}
\begin{algorithm}
\label{alg:reconstruction}
\caption{Cluster graph reconstruction (CGR).}
\SetKwFunction{Compute}{Compute}
\SetKwData{AFalse}{False}
\SetKwFunction{Set}{set}
\SetKwFunction{CI}{CI}
\SetKwFunction{CCT}{CCT}
\SetKwFunction{Stop}{STOP}
\SetKwFunction{Let}{let}
\SetKwFunction{Choose}{choose}
\SetKwInOut{Input}{Input}
\SetKwInOut{Output}{Output}
\Input{Let $U,D, \ve s$ be as in Lemma~\ref{lem:cgr}.}
\Output{Cluster graph $\clustergraph(U,D,S)$.}
\BlankLine
\Let $G_U$ be a Gr\"obner basis of $J_U$ w.r.t. to an arbitrarily chosen term ordering; \\
\Set $N:= \{ \NF(\ve s,G_U) \}$; \\
\Set $V:= \emptyset$; \\
\Set $E:= \emptyset$; \\
\While {$N \neq \emptyset$}{
\Choose $\ve u \in N$; \\
\Set $N := N \setminus \{\ve u\}$; \\
\Set $V := V \cup \{\ve u\}$; \\
\For{$(\ved',\ved'') \in D$}{
	\Set $\ve v:= \CCT(U,\ve u, \ve d')$;\\
	\If{$\ve v \neq$ \AFalse}{
		\If{$\NF( (\ve v-\ve d'+\ved'') - \ve u, G_U) \neq \ve 0$}{
		\Set $\ve w:= \NF(\ve v-\ved'+\ved'', G_U)$; \\
		\Set $E := E \cup \{(\ve u,\ve w)\}$; \\
		\If{$\ve w \not \in V$}{
		\Set $N := N \cup \{\ve w\}$; \\
		\Set $V := V \cup \{\ve w\}$; \\
		}
		}
	}
}
\Return $(V,E)$;
}
\end{algorithm}

We will now explain how the cluster approach can be used to answer Question 1 and Question 2 mentioned in Section~\ref{sec:introduction}. Let $\ves\in\Z^n_+$ and $M=U\cup D\subseteq \Z^n_+\times\Z^n_+$ be as above and let $\clustergraph(U,D,S)$ be the cluster graph that, under the assumption that $J_U:\vex^{\ved'}=J_U$ for all $(\ved',\ved'')\in D$, can be computed by Algorithm~\ref{alg:reconstruction}.

\begin{remark}[Answer to Question 1]
  We can answer Question 1, i.e., given two states $\ve s,\vet\in\Z^n_+$, decide whether $\ves\rightarrow_M \vet$. It suffices to calculate, e.g., the shortest path from $\cluster(\ves)$ to $\cluster(t)$ in $\clustergraph(U,D,S)$. If such a path exists, then we have $\ves\rightarrow_M \vet$. Otherwise such a path does not exist.
\end{remark}

As the cluster graph is usually considerably smaller than the state graph, this test can be performed even for larger reaction networks. If one is interested in a specific path from $\ve s$ to $\vet$ which, e.g., minimizes the energy needed for the reaction, we can attach the required energy as weights to the edges and then use the shortest path algorithm to compute the desired path. In order to answer Question 2, we have to define exactly when we consider two paths {\em essentially different}:

\begin{definition}
  Let $\ves,\vet\in\Z^n_+$ be two states such that $\ves\rightarrow_M \vet$. Further let $(\vev_i)_{i\in [n]}, (\vew_i)_{i\in [m]}\subseteq\Z^n_+$ be two (directed) paths in the state graph $\Gamma_M$ with $\vev_1=\ve w_1=\ve s$ and $\ve v_n=\ve w_m=\ve t$. Then $(\ve v_i)_{i\in [n]}, (\ve w_i)_{i\in [m]}$ are {\em essentially different} if the induced paths in $\clustergraph(U,D,S)$ are different.
\end{definition}

With the definition as above, Question 2 can be answered as follows.

\begin{remark}[Answer to Question 2]
  Let $\ve s, \ve t\in\Z^n_+$ be two states such that $\ves\rightarrow_M \ve t$. Enumerate all paths from $\cluster(\ves)$ to $\cluster(t)$ in $\clustergraph(U,D,S)$ by performing a depth-first search on $\clustergraph(U,D,S)$.
\end{remark}

Usually, enumerating \emph{all} paths connecting two nodes in a graph is rather expensive (actually already path counting is \#P-complete). The reduction to the cluster graph $\clustergraph(U,D,S)$ reduces the problem size though, so that realistic networks can be tackled.

\section{Computational results}
\label{sec:compRes}
We will now provide some computational results for Algorithm~\ref{alg:reconstruction}. All computations were performed in {\mbox CoCoA 4.7} (see \cite{CocoaSystem}) on a machine with a dual core x86\_64 processor with 2Ghz and 2GB of main memory. For the permanganate/oxalic acid reaction and the corresponding mass/charge balance equations (\ref{eq:oxalic}) we obtain $1022$ elementary reactions $M=U\cup D$. We had to construct $19\cdot 18 + 1=343$ Gr\"obner bases for the cluster graph reconstruction. All the necessary Gr\"obner bases computations for this system were performed in less than $2$ min ($1$m$53$s) and the remaining time for the actual cluster graph reconstruction can be neglected (<$1$s) as only normal form computations were performed. The Gr\"obner basis of $J_U$ contained $136$ elements and the Gr\"obner bases of the corresponding ideals $J_U:x_i$ were of similar size.

As shown in \cite{CHEMNETS-1}, $16$ species do suffice to explain the reaction process. We performed the same computations for the reduced matrix and obtained both, significantly reduced computational time and size of the Gr\"obner bases. In this case the $16\cdot 15+1=241$ Gr\"obner bases computations were performed in less than $14$ secs and the Gr\"obner basis of $J_U$ contained $29$ elements. The Gr\"obner bases of the ideals $J_U: x_i$ were again of similar size. Again, the time for computing the cluster graph, see Figure \ref{fig:clustGraph}, is neglectable (<$1$s).

We were also able to compute a Gr\"obner basis of $J_U$ for a large reaction network from an air pollution model involving 80 species and 13556 reversible elementary reactions. The Gr\"obner basis of $J_U$ could be computed in less than $5$ min ($4$m$48$s) and contained $5740$ elements.

\section{Concluding remarks}
\label{sec:conclusion}

We have given algorithms to construct (under certain conditions on $J_U$ and $D$) the induced subgraph $\clustergraph(U,D,S)$ of $\clustergraph(U,D)$ that is reachable from a given finite set $S$ of states. From $\clustergraph(U,D,S)$ we can extract useful information about the decomposition of the overall reaction. For example, transitions can be identified that have to occur in \emph{any} decomposition. Our computations for the permanganate/oxalic acid reaction and the air pollution problem indicate that our algorithms are applicable in practice even for real-world size problems with a large number of species. Therefore, we are confident that our cluster network approach will turn out very useful in the study of more complex biochemical reaction systems.

From a mathematical point of few, some open questions remain:
\begin{itemize}
  \item Find algorithmically all clusters that can be reached from cluster $\cluster(\vey)$ via a given transition $(\ved',\ved'')$. This is an essential step in order to reconstruct the cluster (sub)graph $\clustergraph(U,D,S)$ in the general case, that is, without our assumptions on $J_U$ and $D$.
  \item The reachability problem ``$\ves\rightarrow_M\vet$?'' is generally hard to decide. In the chemical setting above, the situation may be easier. Thus, is it possible to decide ``$\ves\rightarrow_M\vet$?'' in polynomial time for elementary chemical reactions?
  \item Another way to cluster states in $\Z^n_+$ is via strongly
    connected components in the state graph $\Gamma_M$. Again, this
    clustering defines an equivalence relation. This relation can be
    encoded into a certain binomial ideal $J_M$ as before. This ideal,
    however, is given only implicitly via $M=U\cup D$ and so the
    question arises how we can construct generators of $J_M$.
\end{itemize}

\bibliographystyle{plain}
\bibliography{notesBib}

\end{document}